\documentclass[english]{amsart}
\usepackage[T1]{fontenc}
\usepackage{textcomp}
\usepackage{subfig}
\usepackage{amsmath}
\usepackage{amssymb}
\usepackage{amsthm}
\usepackage[english]{babel}
\usepackage{graphicx}
\theoremstyle{definition}

\newtheorem{theorem}{Theorem}
\newtheorem{proposition}[theorem]{Proposition}
\newtheorem{corollary}[theorem]{Corollary}

\newenvironment{proofof}[1]{\noindent {\bf{Proof of #1.}}}{ \hfill\qed\\ }

\def\Tt{T_{\theta}}
\def\A{\mathcal{A}}
\def\L{\mathcal{L}}
\def\BL{\mathcal{BL}}

\title{Polygonal Billiards with One Sided Scattering}
\author {Alexandra Skripchenko}
\address{CNRS UMR 7586, Institut de Math\'ematiques de Jussieu - Paris Rive Gauche, Batiment Sophie Germaine, Case 7021, 75205 Paris Cedex 13, France\\
and Laboratory of Geometric Methods in Mathematical Physics, Lomonosov Moscow State University, Moscow 119991, Russia}
\email{sashaskrip@gmail.com}

\author{Serge Troubetzkoy}
\address{Aix Marseille Universit\'e, CNRS, Centrale Marseille, I2M, UMR
  7373, 13453 Marseille, France}

\curraddr{ I2M, Luminy,\\ Case 907,\\ F-13288 Marseille CEDEX 9,\\ France}

  \email{serge.troubetzkoy@univ-amu.fr}

\begin{document}
\maketitle
\begin{abstract}
We study the billiard on a square billiard table with a one-sided vertical mirror. We associate trajectories of these billiards with double rotations and study orbit behavior and questions of complexity.
 \end{abstract}

\section{Introduction}
The table $\Pi$ we consider consists of the square $[0,1/2]\times [0,1/2]$ with a vertical wall $I$ connecting the points $(a,0)$ and $(a,b)$. We play billiard on this table, with $I$ acting as a one-sided mirror. That is, we consider a point particle in $\Pi$ and a direction $\theta$ in $\mathbb S^{1}$, the particle  moves at unit speed in the direction until it reaches the boundary of the table.
\begin{figure}[h]
\vspace{-3cm}
\includegraphics[width=9cm,height=11.5cm]{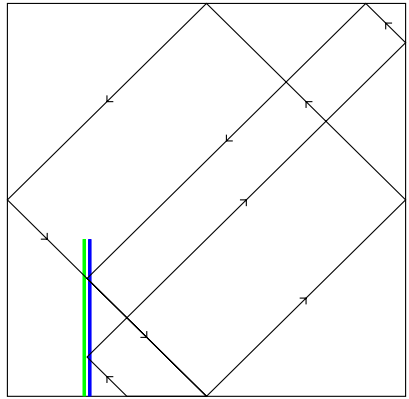}
\vspace{-3.cm}
\caption{The green side is  transparent  while the blue side is reflective.}
\label{1}
\end{figure} 
If it arrives at the left side of $I$ it passes through it unperturbed, while if it arrives at  the right side of $I$ or at the boundary of the square it is reflected with the usual law of geometric optics, the angle of incidence equals the angle of reflection (see Figure \ref{1}). 
Polygonal billiards with one-side straight mirrors were briefly described by M. Boshernitzan and I. Kornfeld \cite{BK} in connection with a special kind of piecewise linear mapping of a semi-interval, called interval translation mappings (ITMs). Interval translation mappings are a natural generalization of interval exchange transformations. 

In this article we prove that billiard flow on such table can be associated with special  interval translation mappings called double rotations (Proposition 1). 
The term  double rotations was introduced in \cite{SIA}, they have also been studied in 
\cite{BT}, \cite{B}, \cite {3}. We show that up to a natural involution there exists a bijection between double  rotations and billiard map we work with, and therefore almost all of our billiard map are of finite type due to the corresponding results for double rotations (Theorem 2 part 2)) (see \cite{SIA} and \cite{3}).
In the other parts of Theorem 2 we collect various interesting implications of this result on unique ergodicity, non-unique ergodicity, and the Hausdorff dimension of the attractor (concretizing a suggestion of Boshernitzan and Kornfeld \cite{BK}).
Our main result is an exact linear formula for complexity of billiard trajectories in a given direction in the case $a=\frac{1}{4}$ (Theorem 4), we also give a linear estimate in case of other rational values of $a$ (Theorem 9). 
The main result also generalizes to certain other rational polygons with one sided scatterers located at an axis of symmetry (Theorem 7). The proof is based on an extension of combinatorial arguments introduced by J. Cassaigne (\cite{C}) for languages with bispecial words (Theorem 11) and also uses a certain symmetry of orbits of the unfolded billiard. 

\section{The results}

A \emph{double rotation} is a map $T: [0,1) \to [0,1)$ of the form:
$$Ty = \left\{
\begin{array}{rl}
y+\alpha  \pmod 1 &\text{if } y \in [0,c)\\
y+\beta \pmod 1 &\mbox{if } y \in [c,1).
\end{array}
\right.$$

\begin{figure}[h]
\vspace{-1.7cm}
\includegraphics[width=9cm,height=11.5cm]{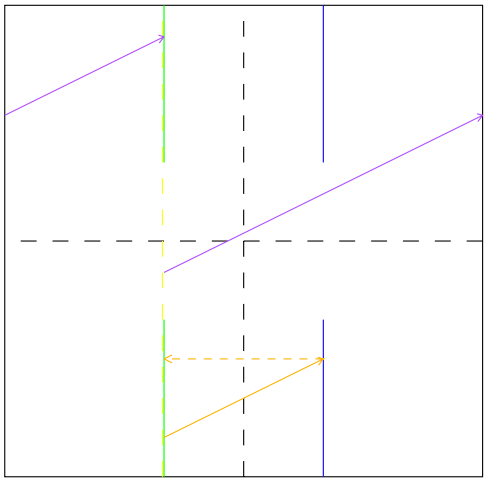}
\vspace{-3cm}
\caption{First return map $T_{\theta}$ is a double rotation}
\label{2}
\end{figure}

Consider the billiard in the table $\Pi$ described in the Introduction. There is a well known construction of unfolding
the billiard in a rational polygon to a translation surface (see for example \cite{MT}).  The same construction applied
to our setting yields  a unit torus consisting
of four copies of the billiard table $\Pi$ with slits corresponding to $I$ which are identified according to the billiard flow:
when we hit a right copy of $I$ (depicted in blue) we jump to the corresponding left copy of $I$ (depicted in green), while the left (green) copies of $I$ are transparent.  With these identifications our billiard is equivalent to the linear flow on the unit slitted torus (see Figure \ref{2}).

Consider the first return map to the vertical section $x=a$ in the unfolded slit torus (see Figure \ref{2}, where this section is drawn by yellow dashed line). We denote this map by $T_{(a,b,\theta)}: \{(x,y): x=a\} \to \{(x,y): x=a\}$. We identify its
domain of definition with the circle $[0,1)$.
It is easy to check that there are only two possibilities for our trajectories (here we assume that the the trajectory is not vertical: $\theta \not \in \{ \pm \frac \pi2 \}$): either the orbit of the point hits the blue wall,  and immediately goes back to the section (this
is depicted by the orange line in Figure \ref{2}); or the orbit of the point does not touch the blue wall (such an orbit is
depicted by the purple line in Figure  \ref{2}). In the first case $T_{(a,b,\theta)}(y) =(1-2a)\tan\theta +y \pmod 1$; in the second case $T_{(a,b,\theta)}(y)=\tan\theta + y \pmod 1$.
Therefore, we have the following:
\begin{proposition}
First return map on the vertical section $x=a$ is a double rotation with the following values of parameters:
\begin{eqnarray*}
\alpha & = & (1-2a)\tan\theta \pmod 1,\\
\beta & = & \tan\theta \pmod 1,\\
c & = & b.
\end{eqnarray*}
\end{proposition}

Let us consider the cubic polynomial $P(x) = x^{3}-x^{2}-3x+1$ and let $\gamma$ be the unique root of $P(x)$ in $[0,1].$
Let us consider the following values of parameters of our billiard: 
\begin{eqnarray*}
a & = & \frac{(1-\gamma)}{2},\\
b & =& \gamma\\
\theta & =& \arctan\gamma.
\end{eqnarray*}
The map $T_{(a,b,\theta)}$ is then the double rotation introduced by Boshernitzan and Kornfeld in \cite{BK}, the
first example of an ITM whose dynamics differs from that of an interval exchange. To explain this difference, and
to collect the most interesting implications of known results to our situation,  we need 
to introduce some notation.
For any ${(\alpha,\beta,c)}$  the {\em attractor} of $T$ is the set $$\Omega :=J\cap T J\cap T^{2} J\cap\cdots,$$ where $J=[0,1).$
If there exists $n\in N$ such that $J\cap T J\cap T^{2} J\cap\cdots\cap T^{n} J=J\cap T J\cap T^{2} J\cap\cdots\cap T^{n+1} J$ then we say that $T$ is of \emph{finite type}.
Otherwise $T$ is of \emph{infinite type}.
Informally, the interval translation map being of finite type means that it can be reduced to interval exchange transformations; if it is of infinite type, then the attractor is a Cantor set. 
The Boshernitzan--Kornfeld example given above is of infinite type.

Fix the parameters $a,b$ and  a direction $\theta$. Let $T:= T_{(a,b,\theta)}$ be the associated double rotation. Let $X_0 := [0,c)$ 
and $X_1 := [c,1)$ be the two intervals of continuity of 
$T$.
The {\em code} $\omega(y) \in \{0,1\}^{\mathbb{N}}$ of the $T$  orbit of $y$ is the sequence of intervals it hits, 
i.e.\ $w(y)_n = i$ if and only if $T^ny$ is in the interval $X_i$.
The {\em language}
 $\mathcal{L} := \mathcal{L}^{a,b,\theta}$ is the set of all infinite codes obtained as $y$ varies, and $\mathcal{L}(n)$ 
 is the set of different {\em words of length $n$}  which appear in
$\mathcal{L}$.  Let $p(n) := p^{a,b,\theta}(n) := \# \mathcal{L}^{a,b,\theta}(n)$.
Note that one could also consider $p^{\infty}(n) \le p(n)$ the number of different words of length
$n$ which appear in $\mathcal{L}$ as $x$ varies in the attractor.

We must exclude directions for which there is a billiard orbit from an end point of $I$
to an end point of $I$.  We call such directions {\em exceptional}. There are
at most a countably many exceptional  directions  since for any positive constant $N$, there is a finite number of 
billiard orbits which start and end at end points of $I$ and have length at most $N$.

The following theorem follows directly from several known results about double rotations, interval translation mappings and piecewise continious interval maps. 
 
\begin{theorem}{\parbox[t]{0.2\linewidth}{$\phantom{1}$  }}\label{thm1}

\noindent
1) For all $(a,b,\theta)$ with $\theta$ non-exceptional the billiard/double rotation is minimal.\\
2)  For almost all  $(a,b,\theta)$ the double rotation is of finite type.\\
3) There exists an uncountable set of $(a,b,\theta)$ so that the Hausdorff dimension of the closure of the attractor is zero, in particular the double rotation is of infinite type.\\
4)  For each $(a,b,\theta)$ with $\theta$ non-exceptional the  billiard/double rotation has at most two ergodic invariant measures.\\
5) There exists an uncountable set of $(a,b,\theta)$ with $\theta$ non-exceptional such  that the billiard/double rotation is not uniquely ergodic.\\
6) For all $(a,b,\theta)$ with $\theta$ non-exceptional the complexity $p(n)$ grows at most polynomially with degree $3$.
\end{theorem}

In \cite{BK} it was suggested that there may exist configurations on a rational billiard table with mirrors which force the light to get concentrated in some arbitrarily small portions of the table. This suggestion is confirmed by
Part 3) of the theorem.

\begin{proof}
It is easy to check that varying the parameters $(a,b,\theta)$ we can obtain all double rotations with only one restriction: 
$\alpha\le\beta.$ There exists an involution between two parts of parameter space ($\alpha\le\beta$ and $\beta\le\alpha$) and the orbit behavior of ITM from these two parts are completely the same: the involution is measure preserving and does not change the dynamics.  Thus all known results on double
rotations hold in our setting. 
1)  We think of a double rotation as  an ITM defined on an interval. First suppose that  this interval is $[0,1)$  (in the coordinates of the definition of double rotation). Either this is an ITM on with 3 intervals of continuity, or 
it  has 4 intervals of continuity. In the later case 
 we choose the origin to be the point $c$, and the double rotation now 
has at most 3 intervals of continuity, thus we can always choose coordinates so that the ITM has at most 3 intervals and
part 1)  follows from Theorem 2.4 of \cite{ST}.  Part
2)  follows from \cite{SIA} Theorem 4.1 or from Theorem 1 of \cite{3},
part 3)  follows from Theorem 10 of \cite{BT},
part 4) follows from Theorem 3 of \cite{BH},
part 5)  follows from Theorem 11 of \cite{BT}, and 
part 6)  follows Theorem 1 of \cite{Ba} (see also Corollary 8 of \cite{BH}).
 \end{proof}

The \emph{billiard flow} $\phi_{t}$ is defined on the phase space $\tilde{\Pi} := \Pi\times S^{1}$, with proper identifications on the boundary, and the \emph{attractor} of the billiard flow is 
$$\bigcap_{t\geq 0}\phi_{t}{\tilde{\Pi}}.$$
An immediate application of the Fubini Theorem to Theorem 2 part 2) yields

\begin{corollary}
For almost every $a,b$ the billiard attractor has full measure.
\end{corollary}

\noindent {\bf Question:} For which polygons with one sided mirrors does the billiard attractor have full measure? Positive measure? Zero measure? 

The main results of this article are improvements of part 6) of the theorem.
We begin with the case  $a =\frac 1 4$, i.e.\ the one sided mirror is in the middle of the square, where we get
a complete description of the complexity.
We define three sets of directions, in the below definitions all the parameters $k_i$ and $l_i$ are integers.
\begin{eqnarray*}
\tilde{A}_1 := \Big \{\theta: \exists k_1, k_2, l_1, l_2:  \forall  \tilde{k_1}\leq k_1,\tilde{k_2}\leq k_2, \tilde{l_1}\leq l_1, \tilde{l_2}\leq l_2:\\
\tan \theta \in  \left (\frac {2l_1}{1+2k_1},  \frac {2l_1+4b}{2k_1+1}\right ) \cap \bigcap_{\tilde{k_1},\tilde{l_1}} \left(\frac {\tilde{l_1}+2b}{\tilde{k_1}+1},\frac {\tilde{l_1}+1}{\tilde{k_1}+1}\right) \cap \\
  \left (\frac {2l_2+2-4b}{1+2k_2},  \frac {2l_2+2}{2k_2+1}\right ) \cap \bigcap_{\tilde{k_2},\tilde{l_2}}\left(\frac {\tilde{l_2}+1}{\tilde{k_2}+1},\frac {\tilde{l_2}+2-2b}{\tilde{k_2}+1}\right)\Big\},
\end{eqnarray*}
\begin{eqnarray*}
\tilde{A}_2 :=  \Big \{\theta: \exists k_1, k_2, l_1, l_2:  \forall  \tilde{k_1}\leq k_1,\tilde{k_2}\leq k_2, \tilde{l_1}\leq l_1, \tilde{l_2}\leq l_2:\\
\tan \theta \in  \left (\frac {2l_1}{1+2k_1},  \frac {2l_1+4b}{2k_1+1}\right ) \cap \bigcap_{\tilde{k_1},\tilde{l_1}} \left(\frac {\tilde{l_1}+2b}{\tilde{k_1}+1},\frac {\tilde{l_1}+1}{\tilde{k_1}+1}\right) \cap \\
 \left (\frac {1+l_2-2b}{1+k_2},  \frac {1+l_2}{k_2+1}\right ) \cap \bigcap_{\tilde{k_2},\tilde{l_2}}\left(\frac {2\tilde{l_2}}{\tilde{k_2}+1},\frac {\tilde{l_2}+2-4b}{\tilde{k_2}+1}\right)\} \bigcup \\
 \Big \{\theta: \exists k'_1, k'_2, l'_1, l'_2:  \forall  \tilde{k'_1}\leq k'_1,\tilde{k'_2}\leq k'_2, \tilde{l'_1}\leq l'_1, \tilde{l'_2}\leq l'_2:\\
 \tan \theta \in  \left (\frac {l'_1+1}{1+k'_1},  \frac {1+l'_1+2b}{k'_1+1}\right ) \cap \bigcap_{\tilde{k'_1},\tilde{l'_1}} \left(\frac {2\tilde{l'_1}+2b}{2\tilde{k'_1}+1},\frac {2\tilde{l'_1}+2}{2\tilde{k'_1}+1}\right) \cap \\
 \left (\frac {2l'_2+2-4b}{1+2k'_2},  \frac {2l'_2+2}{2k'_2+1}\right ) \cap \bigcap_{\tilde{k'_2},\tilde{l'_2}}\left(\frac {\tilde{l'_2}+1}{\tilde{k'_2}+1},\frac {\tilde{l'_2}+2-2b}{\tilde{k'_2}+1}\right)\Big\}   
\end{eqnarray*}
\begin{eqnarray*}
\tilde{A}_3 :=  \Big \{\theta: \exists k_1, k_2, l_1, l_2:  \forall  \tilde{k_1}\leq k_1,\tilde{k_2}\leq k_2, \tilde{l_1}\leq l_1, \tilde{l_2}\leq l_2:\\
\tan \theta \in  \left (\frac {l_1+1}{1+k_1},  \frac {1+l_1+2b}{k_1+1}\right ) \cap \bigcap_{\tilde{k_1},\tilde{l_1}} \left(\frac {2\tilde{l_1}+2b}{2\tilde{k_1}+1},\frac {2\tilde{l_1}+2}{2\tilde{k_1}+1}\right) \cap \\
  \left (\frac {1+l_2-2b}{1+k_2},  \frac {1+l_2}{k_2+1}\right ) \cap \bigcap_{\tilde{k_2},\tilde{l_2}}\left(\frac {2\tilde{l_2}}{\tilde{k_2}+1},\frac {\tilde{l_2}+2-4b}{\tilde{k_2}+1}\right)\Big\},   
\end{eqnarray*}

For $i= 1,2,3$ let $A_i := \tilde{A}_i \setminus   \{\theta \text{ exceptional}\}$.
The set are pairwise disjoint and 
$$\bigcup_{i=1}^3 A_i = \mathbb{S}^1 \setminus   \{\theta \text{ exceptional}\}.$$

\noindent We give a geometric description of these sets. We work in the unfolded model   (Figure \ref{2}) and consider the two singular points of the forward billiard map. We consider the backward trajectories of these two singular points, the sets are distinguished by the fact which part of the boundary of the semi-transparent mirror is first touched by the backwards trajectory, the transparent side or the reflecting side. This identifies how many preimages a corresponding point of the orbit has (none or two). In particular, $A_{3}$ corresponds to the case when both of the trajectories first touched the reflecting side, $A_{1}$ to the case when both of the trajectories first touched the transparent side and $A_{2}$ is the intermediate case, one trajectory first touches the reflecting side and the other tranjectory first touches the transparent side. 

\begin{theorem}\label{thm2}
Suppose $a =\frac 1 4$ and that
$\theta$ is non-exceptional. Then either, $\theta \in  A_1$, 
and $$p(n) = n+1 \quad \text{for all} \ n \ge 0,$$ 
or  there exists a positive
integer constant $C_{\theta}$ so that
\begin{enumerate}
\item for $\theta \in A_2$, $p(n) = 2n - C_{\theta}$ for all sufficiently large $n$, or
\item for $\theta \in A_3$,  $p(n) = 3n - C_{\theta}$ for all sufficiently large $n$.
\end{enumerate}
\end{theorem}

Remark, in cases 1 and 2  the double rotation $\Tt$ is of infinite type where by $\Tt$ we mean $T_{1/4,b,\theta}$ for any arbitrary $b$.
The behavior of the complexity for small $n$ will also be described in the proof.
Since the map $\Tt$ is minimal, we can apply a theorem of Boshernitzan, \emph{a minimal symbolic system with $\lim_{n\to\infty} sup \frac{p(n)}{n}<3$ is uniquely ergodic} (\cite{FM}[Theorem 7.3.3]) to conclude that
\begin{corollary}
If $\theta \in A_1 \cup A_2$ then $\Tt$ is uniquely ergodic. 
\end{corollary}

Note that $\tilde{A}_2$ is open, and exceptional directions are countable, thus we have shown that
for all but countably many $\theta$ in an open
set of $\theta$ the map  $\Tt$ is of infinite type, uniquely ergodic and of linear complexity.

\begin{corollary}
In the case $a=\frac{1}{4},$ the billiard attractor has positive, but not full measure.
\end{corollary}

Except for the exact computation of the sets $A_i$, Theorem \ref{thm2} is a special case of a
more general result.  Take any rational polygon,  reflect it in one side.  Erase part of the side, make the
other part a one sided mirror (which will be part of the line $x=0$ to be concrete), to produce a table $P$. Consider the slitted flat surface $M := M(P)$ associated
with $P$ (see for example \cite{MT}) with the slits identified as in the square case.
For the moment consider $P$ without the one-sided mirror, and $M$ without the slits, let 
 $\frac{m_i}{n_i} \pi$ be the angles between the sides of $P$, and let $N$ be the least common multiple of the
$n_i$.  Then a standard computation
shows that $M(P)$ has $R := N \sum \frac 1 {n_i}$ vertices (not counting the endpoints of the slits) \cite{MT}.

Consider  the
section $x=0$ as a subset of $P$.
In $M$ there are $2N$ copies of this section.  If we fix $\theta$ non-exceptional (no orbit from any vertex of $P$ to any vertex of
$P$) then there are $N$ copies of the section for which
the linear flow on $M$ jumps via an identification (the linear flow goes through the other $N$ copies of the
section as if they where not there). We consider the first return map $\Tt$ to the $N$ copies which produce a 
jump.  The top and bottom of the section have already been counted as vertices of $P$.
The map $\Tt$ is an interval translation map with $R+N$  intervals of continuity, the $R$ coming from the $R$ singular
points on $M$, and the  $N$ from the end point of the
one-sided mirror (which we assume starts at the bottom of $P$).  

We assign a symbol to each of these intervals and code the billiard orbit by these $R$ symbols.
Let $p(n)$ be the complexity of the associated language.  

\begin{theorem}\label{thm3}
If $P$ as above and $\theta$ is non-exceptional then there exists a positive constant
$C_{\theta}$ and $k \in \{0,1,2,\cdots,R+N\}$ so that
$p(n) = (R+ N-1 + k)n - C_{\theta}$ for all sufficiently large $n$.  
\end{theorem}

\begin{figure}[h]
\centering
\vspace{-3cm}
\hspace{-1cm} \mbox{\subfloat{
\includegraphics[scale=0.3]{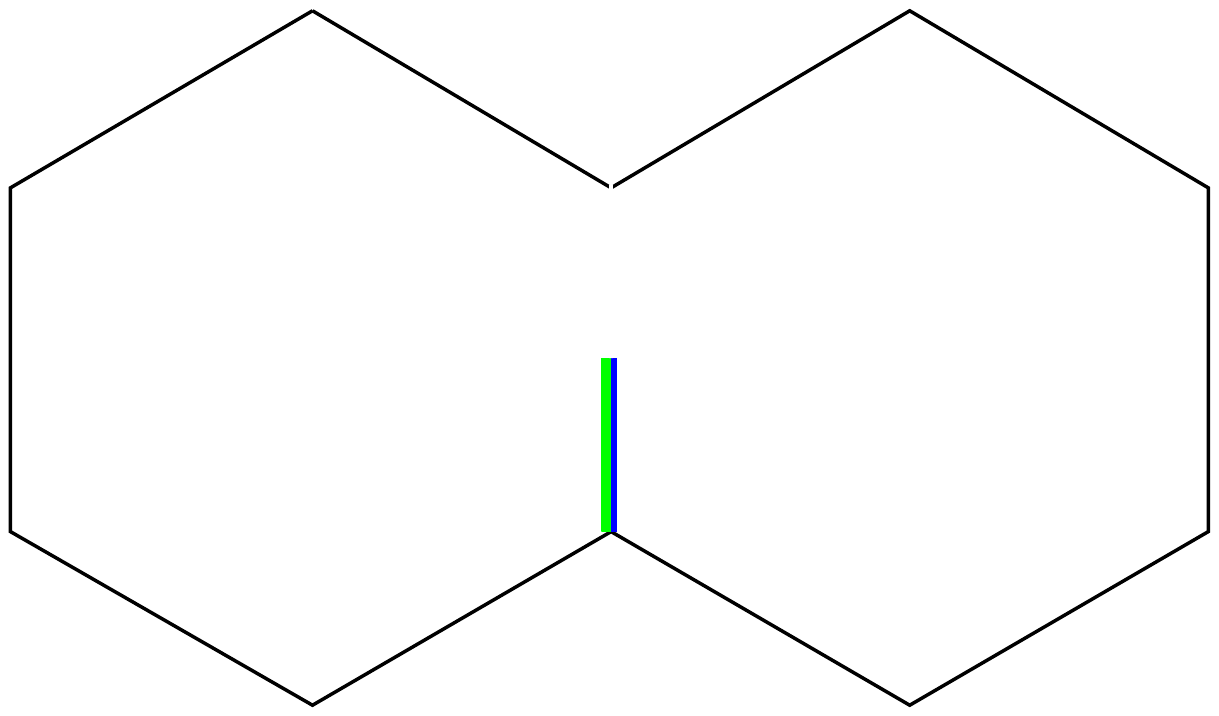}}
\quad
\subfloat{
\includegraphics[scale=0.3]{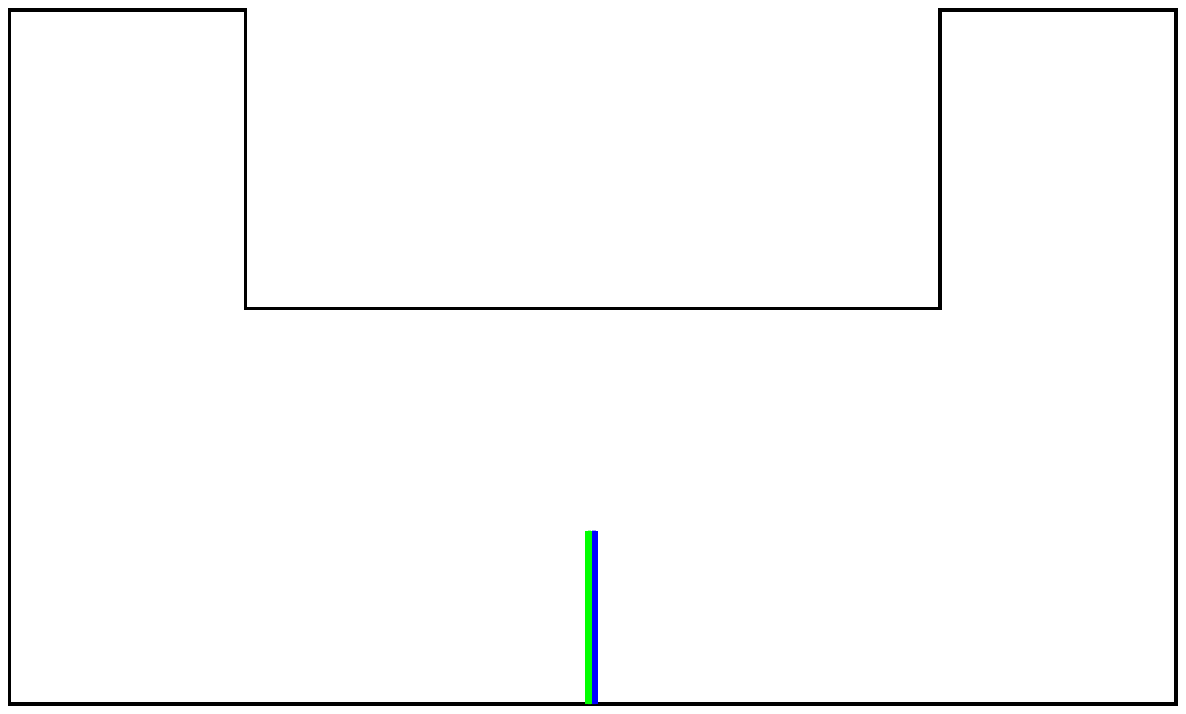}}}
\vspace{-3cm}
\caption{Polygons with  linear complexity. }
\label{fig3}
\end{figure}
Example 1:
In the double hexagon (Figure \ref{fig3}) all angles have $n_{i}=3$, so $N=3, R=18$. In this case all $m_i \ne 1$, 
so we do not have any removable singularities.

Example 2:
In the $U$ shaped figure (Figure \ref{fig3}) one can calculate $N=2, R=12$.  However, on $M$ only 4 of these singularities are
actual singularities (the two points with angle $\frac 3 2 \pi$ and the two copies of the end point of the one sided mirror
on $M$), the others are removable singularities. We also remark that in the original case of the square $N=2, R=4$ and there are 4 removable singularities.

We remind the reader of two results: \emph{an aperiodic ITM with r intervals of continuity has at most $\{r/2\}$ minimal components} (Theorem 2.4 in \cite{ST}) and \emph{if $K\geq 3$ is an integer, then a minimal symbolic system satisfying} $lim sup_{n\to\infty}\frac{p(n)}{n}<K$ \emph{admits at most K-2 ergodic invariant measures} (Theorem 7.3.4 in \cite{FM}).
Combining these two results with Theorem 7 yields
\begin{corollary}
There are at most $[(R+N)/2]$ minimal components and 
on each minimal component the number of ergodic invariant measures for the ITM/billiard is at most $2N+2R-2$.
\end{corollary}

We turn to the question of complexity for rational position of the mirror in the square case, i.e.\ $a \in \mathbb{Q}$.
Suppose $a =\frac p q$.  
We consider the
return map to the  vertical sections $x= \frac k q$ with $k \in \{0, 1, \dots q-1\}$. The first return map $S$ to these vertical sections is  given by
\
$$S(i,y) = \left\{
\begin{array}{ll}
(p,y+\tan{\theta}/q \hspace{-.2cm} \pmod 1) & \mbox{if } i = q-p \text{ and } y \in [-b,b],\\
(i+1   \hspace{-.2cm}  \pmod q,y+\tan{\theta}/q  \hspace{-.2cm}  \pmod 1)  &\mbox{otherwise}.
\end{array}
\right.$$
We code orbits of this map by a 3 letter alphabet, 2 letters for the section $x = a$ (where the map is discontinuous)
 and a third letter for all
the other sections.

\begin{theorem}\label{thm5}
If $a =\frac p q$, $\theta$ is non-exceptional then $p_{\theta}(n) \le (2+2q) n$ for all $n$.
\end{theorem}

Theorem 2.4 in \cite{ST} and Theorem 7.3.4 in \cite{FM} imply that 
\begin{corollary}
There are at most $q$ minimal components and 
on each minimal component the number of ergodic invariant measures for the ITM/billiard is at most $1+2q$.
\end{corollary}

\section{Cassaigne's formula}

The main technical tool will be a variant of Cassaigne's formula \cite{C}.
Let $\A$ be a finite alphabet, $\L \subset \A^{\mathbb{N}}$ be a language, $\L(n)$ the set of words of length $n$ which appear in $\L$, and $p(n) := \# \L(n)$.
Note that $p(0) = \# \{\emptyset\} = 1$.
For any $n \ge 0$ let $s(n) := p(n+1) - p(n),$ and thus $$p(n) = 1 + \sum_{i=0}^{n-1} s(i).$$
For $u \in \L(n)$ let
\begin{eqnarray*}
m_l(u) & := & \#\{a \in \A: au \in \L(n+1)\},\\
m_r(u) & := & \#\{b \in \A: ub \in \L(n+1)\},\\
m_b(u) & := & \#\{(a,b) \in \A^2: aub \in \L(n+2)\}.
\end{eqnarray*}
We remark that while $m_r(u) \ge 1$ the other two quantities can be $0$.
A word $u \in \L(n)$ is called left special if $m_l(u) > 1$, right special
if $m_r(u) > 1$ and bispecial if it is left and right special.
Let $\BL(n) := \{u \in \L(n): u \text{ is bispecial}\}.$
Let $\L_{np}(n):= \{v \in \L(n): m_l(v) = 0 \}$.

In this section we show that 

\begin{theorem}\label{thm::com}
$$s(n+1) - s(n) = 
\sum_{v \in \BL(n)} \Big (m_b(v) - m_l(v) - m_r(v) + 1 \Big )   - \sum_{v \in \L_{np}(n):m_r(v) >1}  
\Big (m_r(v) - 1 \Big ). $$
\label{combi}
\end{theorem}

{\bf Remark:} Cassaigne proved this theorem in the case of recurrent languages (i.e. $\L_{np}(n) = \emptyset$)  \cite{C} (see
\cite{CHT} for a English version of the proof).  We use the same strategy of proof.

\begin{proof}
Since for every $u \in \L(n+1)$ there exist $b \in \A$ and
$v \in \L(n)$ such that $u = vb$ we have
$$s(n) = \sum_{u \in \L(n)} (m_r(u) - 1).$$
Thus
$$s(n + 1) - s(n)  = \sum_{u \in \L(n+1)} \Big (m_r(u) - 1 \Big ) - 
\sum_{v \in \L(n)} \Big (m_r(v) - 1 \Big ).$$
Now, for $u \in \L(n+1)$ we can write $u = av$ with $a \in \A$ and $v \in \L(n)$. 
Let $\L_p(n):= \{v \in \L(n): m_l(v) \ge 1\}$.
Thus
$$s(n+1) = \sum_{v \in \L_p(n)} \left [ \sum_{av \in \L(n+1)}
\Big (m_r(av) - 1 \Big ) \right ].$$
Let $\L_{np}(n) := \L(n) \setminus \L_p(n)$, and thus $s(n+1) - s(n)$ equals
$$ \sum_{v \in \L_p(n)} \left [ \sum_{av \in \L(n+1)}
\Big (m_r(av) - 1 \Big )    -    \Big ( m_r(v) - 1 \Big )    \right ]  - \sum_{v \in \L_{np}(n)}  
\Big (m_r(v) - 1 \Big ). $$

 For any word $v \in \L_p(n)$ with $av \in \L(n+1)$ any legal prolongation
to the right of $av$ is a legal prolongation to the right of $v$ as well
thus if $m_r(v) = 1$ then $m_r(av) = 1$.  Thus words with $m_r(v) = 1$
do not contribute to any of the above sums.
Thus $s(n+1) - s(n)$ is equal to the above sum
restricted to those $v$ such that $m_r(v) > 1.$ 
For the left sum, if furthermore $m_l(v) = 1$ then there is only a single $a$ such
that $av \in \L(n+1)$. For this $a$ we have $m_r(av) = m_r(v)$ thus
such words do not contribute to the left sum. Thus the only terms
which contribute to the left sum are the bispecial words, and to the right  the
words for which $m_r(v) > 1$; in other words $s(n+1) - s(n)$ equals
$$ \sum_{v \in \BL(n)} \left [ \sum_{av \in \L(n+1)}
\Big (m_r(av) - 1 \Big )    -    \Big ( m_r(v) - 1 \Big )    \right ]  - \sum_{v \in \L_{np}(n):m_r(v) >1}  
\Big (m_r(v) - 1 \Big ). $$

For any $v \in \BL(n)$ we have
$$ m_b(v) = \sum_{av \in \L(n+1)} m_r(av)$$
and
$$ m_l(v) = \sum_{av \in \L(n+1)} 1,$$
thus $s(n+1) - s(n)$ equals 

$$ \sum_{v \in \BL(n)} \Big (m_b(v) - m_l(v) - m_r(v) + 1 \Big ) - \sum_{v \in \L_{np}(n):m_r(v) >1}  
\Big (m_r(v) - 1 \Big ). $$
\end{proof}

\section{The proofs}
\begin{proofof}{Theorem \ref{thm2}}
We use Theorem \ref{combi}. 
In our case $m_r(v) = 1$ or $2$, so the second term of the equation reduces to $\#\{  v \in \L_{np}(n):m_r(v) >1 \}$, thus 
$s(n+1) - s(n)$ equals
$$\sum_{v \in \BL(n)} \Big (m_b(v) - m_l(v) - m_r(v) + 1 \Big ) - \#\{  v \in \L_{np}(n):m_r(v) >1 \}. $$

We have $p(0) = 1$ and $p(1) = 2$.  
Suppose that $y$ has two preimages, 
and that $\Tt^n(y)$ is in the boundary of one of the two intervals of continuity of $\Tt$, i.e.\ the billiard orbit of $y$ arrives at the top of the reflecting side of $I$ after $n$ steps. 
Consider the code $v$ of
length $n$. Clearly $m_l(v) = m_r(v) = 2$. Since $\theta$ is non-exceptional then $m_b(v) = 4$.
(Note that if $\theta$ is exceptional then
$m_b(v) = 3$ and thus the orbits would not contribute to the sum.)
 Thus for non-exceptional directions 
 $$s(n+1) - s(n)  =  \# \BL(n) - \#\{  v \in \L_{np}(n):m_r(v) >1 \}.$$
 
We switch back and forth between the language of double rotations and that of the billiard.
Let $\{e,f\}$ be the  two (common) endpoints of $X_0$ and $X_1$ (in the original definition of double rotations these points are called $0$ and $c$). Let $\{e_{i}^{(n)}\} := \Tt^{-n}(e)$ be the collection of nth-preimages of $e$ which we will denote by $T_e(n)$ for short.  Then we will define the tree $T_e$ of preimages of $e$
to be the set $$T_e = \cup_{n \ge 0}T_e(n),$$ with a directed arrow from $e_{i}^{n}$ to $e_{j}^{n-1}$ if $\Tt e_i^{n} = e_j^{n-1}.$ We define similarly the preimage tree $T_{f}$.

Any right special word  corresponds to a billiard orbit which hits $e$ or $f$, thus we can decompose
$s(n+1) - s(n)$ into two parts, those words contributing to this difference corresponding to a  billiard orbit
arrives at $e$, and those which arrive at $f$; we note these two contributions by
$$(s(n+1) - s(n))_e  =  \# \BL_e(n) - \#\{  v \in \L_{np}(n):m_r(v) >1 \}_e$$
for the point $e$, and similarly for the point $f$. This formula can be seen as counting the number of leaves at level $n$ of a weighted tree, a vertex $e^{n}_{i}$ has weight $+1$ if it has two preimages (contributing one element to  $\BL_e(n)$), weight $0$ if it has one preimage, weight $-1$ if it has no preimages (contributing one element  in $\L_{np}(n)$) (Figure \ref{Tree}).

\begin{figure}[h]
\vspace{-4cm}
\includegraphics[width=9cm,height=11.5cm]{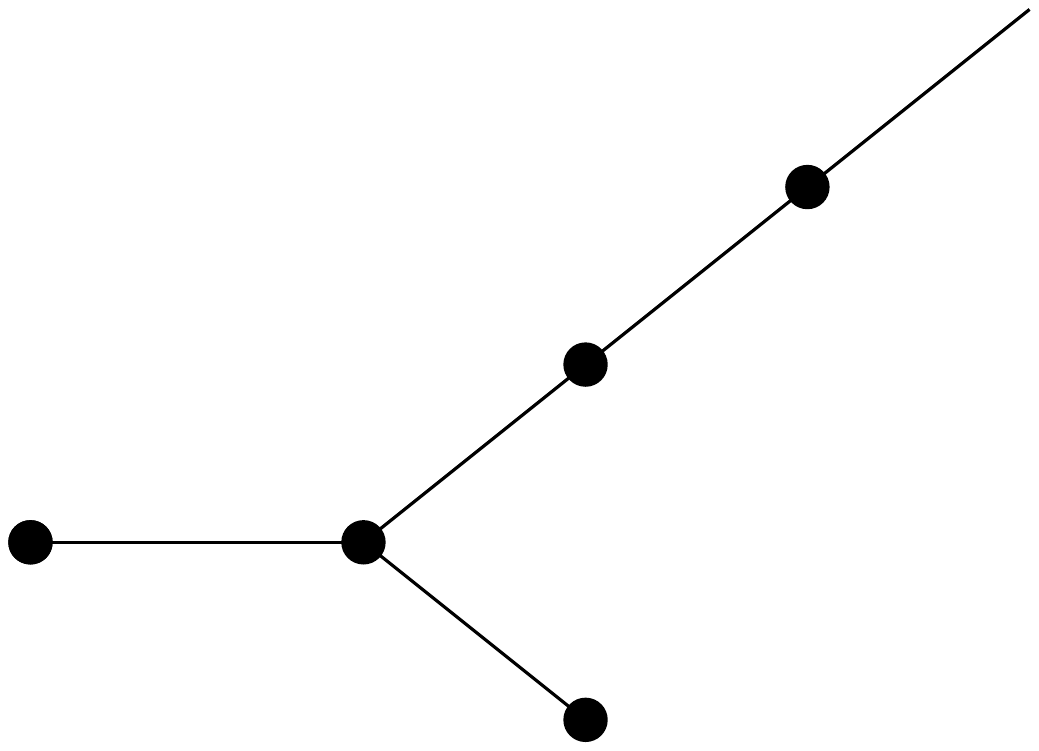}
\put(-152,153){$e$}
\put(-120,153){$+1$}
\put(-91,168){$0$}
\put(-64,187){$0$}
\put(-90,132){$-1$}
\vspace{-4cm}
\caption{The preimage tree of $e$.}
\label{Tree}
\end{figure} 

For any $n$ such that each of the $e^{n}_{j}$ has a unique preimage, we have $\# \BL_e(n) = \#\{v\in\L_{np}: m_{r}(v)>1\}_{e}=0$ and thus $(s(n_0+1) - s(n_0))_e = 1$.

Now consider an $n$ so that at least one point  $e_j^{n-1}$ has two preimages.
 The first time this happens, $n_0 := n_0^e$, if this never happens
then set $n_0 := +\infty$.
Thus $\# \BL_e(n_0) = 1$ and  $ \#\{  v \in \L_{np}(n_0):m_r(v) >1 \}_e=  0$. This implies
$(s(n_0+1) - s(n_0))_e = 1$.
By the symmetry, the next time, $n_1^e$ one of these backwards orbits hits the mirror, the other one will also
hit the mirror from the other side.  Thus one of the orbits will have two preimages and the other will be non-left-extendable.  Thus   $\# \BL_e(n_1^e) = \#\{  v \in \L_{np}(n_1^a):m_r(v) >1 \}_e=  1$.  The same holds for all times
$n_k^e$ for which $\# BL_e(n_k^e) > 0$.  We have thus shown that $n_0^e$ is the unique time for which
$(s(n+1) - s(n))_e \ne 0$. 
The same holds for $(s(n+1) - s(n))_f \ne 0$, but $n_0^e$ and $n_0^f$ are not necessarily equal.

If $n_0^e = n_0^f =:  n_0$ then
\begin{equation*}
s(n+1) - s(n) = 
\begin{cases} 0 & \text{if $n \ne n_0$,}
\\
2 &\text{if $n  = n_0$.}
\end{cases}
\end{equation*}
By definition $p(0) = 1$ and thus $s(0) = p(1) - p(0) = 2 - 1 = 1$.
Then $s(n) = s(0) = 1$ for all $n \le n_0$ and
$s(n) = 3$ for all $n > n_0$.  Thus 
\begin{equation}\label{c1}
p(n) = 
\begin{cases} n+1 & \text{if $n \le n_0$,}
\\
3n+(1-2n_0)  &\text{if $n > n_0$.}
\end{cases}
\end{equation}
In particular if $n_0 = \infty$, i.e.\ if the preimages of both points $e$ and $f$ disappear before being
doubled, then 
$$p(n) = n+1 \qquad \text{for all} \ n \ge 0.$$
It is easy to check that $n_0 = \infty$ happens exactly when $\theta \in A_1$.

On the other hand if $n_0^e \ne n_0^f$,
 if follows that
\begin{equation*}
s(n+1) - s(n) = 
\begin{cases} 0 & \text{if $n \not \in \{n_0^e,n_0^f\}$,}
\\
1 &\text{if $n  \in \{n_0^e,n_0^f\}$.}
\end{cases}
\end{equation*}
Set
$N_0^- := \min(n_0^e,n_0^f)$ and $N_0^+ := \max(n_0^e,n_0^f)$. 
We have $s(n) = s(0) = 1$ for all $n \le N_0 $,
$s(n) =  2$ for all $n \in (N_0^-,N_0^+]$,
$s(n) = 3$ otherwise.  Thus 
\begin{equation}\label{c2}
p(n) = 
\begin{cases} n+1 & \text{if $n \le N_0^-$,}
\\
2n + (1 - N_0^-)  &\text{if $n \in (N_0^-,N_0^+]$,}
\\
3n + (1- N_0^- - N_0^+) & \text{otherwise}.
\end{cases}
\end{equation}
If $N_0^+ = \infty$ then exactly one of the two points $e$ or $f$ disappears before being doubled, thus $\theta\in A_{2},$ the ITM
is of infinite type and
\begin{equation*}
p(n) = 
\begin{cases} n+1 & \text{if $n \le N_0^-$,}
\\
2n + (1 - N_0^-)  & \text{if $n >  N_0^-$.}
\end{cases}
\end{equation*}
If $N_0^+ < \infty$ then neither preimage disappers before being doubled, thus $\theta\in A_{3}$.
\end{proofof}

\begin{proofof}{Theorem \ref{thm3}} The main difference with Theorem \ref{thm2} is that instead of two points $\{e,f\}$
which  produce right special words, there are now $R+N$ points. Call these point $\{g_1,g_2,\dots,g_{R+N}\}$. 
The other difference is that there are $R+N$ intervals of continuity, thus $p(1)=R+N$ and thus $s(0) = p(1)  - p(0) =  R+N - 1$.  

Otherwise the proof is identical, for each of the $g_i$ we construct its preimage tree.
We consider the first time $n_{g_i}$ when the preimage is doubled.  The symmetry argument
is identical to the square case, once a $g_i$ has two preimages at some time,  it has two preimages for
larger times. The $k$ in the statement of the theorem then corresponds to the number of $g_i$ which
have two-preimages at a certain time, while for the other $R+N-k$ points the preimages disappear before being doubled.
As before, $C_{\theta}$ is responsible for the events that happened before the mo-
ment when the first preimage was doubled.
We conclude that $s(n+1) - s(n) =k$ for $n$ sufficiently large and thus
 $p(n) = (s(0) + k) n - C_{\theta}  = (R+N -1 + k)n - C_{\theta}$.
\end{proofof}

\begin{proofof}{Theorem \ref{thm5}}
Remember that in this case the map $S$ is the first return map to the union of vertical sections $x=\frac{k}{q}$ described before the Theorem.
The proof follows the same line as the the case $a = \frac 1 4$, in which we
argued that by symmetry that once there are two preimages of $e$, each time one disappears a new one
appears.  This is no longer true in the general rational case.
As in the case $a = \frac 1 4$ we need to consider the preimage tree under the map $S$ of the points $e$ and $f$ (which
as in case $a=\frac 1 4$ are the points of discontinuity on the circle $x=a=\frac{p}{q}$), the map being continuous on the other circles.

As mentioned in the proof of Theorem 4, the difference $(s_{n}-s_{0})_{e}$ counts the number of leaves of level $n$, but 
since all points in the set $\{S^{-n}(i,y)\}$ have the same second coordinate, we have $\# \{S^{-n}(i,y)\} \le q$, for any $n \ge 0$, for any point $(i,y)$. Thus 
$$(s(n+1)-s(n))_e \leq q$$ for all $n$.

Taking into account the contribution of $e$ and $f$ yields 
$$s(n) - s(0)  \le 2q$$
and since $s(0) = p(1) - p(0) = 3 -1 = 2$,  we conclude that 
 $$p(n) = 1 + \sum_{i=0}^{n-1} s(i) = 3 + \sum_{i=1}^{n-1} s(i)   \le 3 + (2 + 2q) (n-1) < (2+2q)n.$$
\end{proofof}

\begin{proposition}\label{prop1}
If $\{a,b\} \cap \Omega= \emptyset$ then $p^{\infty}(n) = n+1$.
\end{proposition}

\begin{proofof}{Proposition \ref{prop1}}
By assumption there are no right special words, thus no bispecial words.  Thus $s(n) = const$.
We have $s(0) = p(1) - p(0) = 3 - 1 = 2$,  thus $s(n) = 1$ and $p(n+1)= p(n) +1.$
It follows that $p^{\infty}(n) = n+1$.   
\end{proofof}

\section{Acknowledgements.} We thank Pascal Hubert for fruitful discussions and Henk Bruin for a critical reading of an earlier 
version of this manuscript. We  gratefully acknowledge the support of
ANR Perturbations and ERC Starting Grant \textquotedblleft Quasiperiodic\textquotedblright. We also thank the anonymous referee for several useful suggestions which improved our presentation.

\end{document}